\documentclass[11pt]{article}
\usepackage{amsmath,amssymb,color,graphicx,bbm}

\newcommand {\C}{\mathbb{C}}

\newcommand{\diag}{\mathrm{diag}}

\newcommand{\matlab}{{\sc matlab}}

\newcommand{\tr}{\mathrm{tr}}
\newcommand{\Res}{\mathrm{Res}}

\newcommand{\beq}{\begin{equation}}
\newcommand{\eeq}{\end{equation}}

  \newtheorem{theorem}                {Theorem}
  \newtheorem{lemma}      {Lemma}
  
  \newtheorem{assumption}{Assumption}
  \newtheorem{remark}{Remark}

\DeclareMathOperator{\T}{T}

\title{First-order Perturbation Theory for \\Eigenvalues and Eigenvectors}
\author{Anne Greenbaum\thanks{Department of Applied Mathematics, University of Washington.}
\and
Ren-cang Li\thanks{Department of Mathematics, University of Texas at Arlington. Supported in part by
 National Science Foundation Grants CCF-1527104 and DMS-1719620.}
\and 
Michael L.~Overton\thanks{Courant Institute of Mathematical Sciences, New York University. Supported in
part by National Science Foundation Grant DMS-1620083.}
}

\date{\today\\
\ \\
Dedicated to Peter Lancaster and G.W.\ (Pete) Stewart\\
{\normalsize Masters of Analytic Perturbation Theory and Numerical Linear Algebra\\
on the Occasion of their 90th and 79th Birthdays}
}

\begin{document}
\maketitle

\begin{abstract} 
We present first-order perturbation analysis of a simple eigenvalue and the corresponding right and left eigenvectors of a 
general square matrix, not assumed to be Hermitian or normal. The eigenvalue result is well known to a broad scientific community. 
The treatment of eigenvectors is more complicated, with a perturbation theory that is not so well known outside
a community of specialists. We give two different proofs of the main eigenvector perturbation theorem. The first, 
a block-diagonalization technique inspired by the numerical
linear algebra research community and based on the implicit function theorem, 
has apparently not appeared in the literature in this form.
The second, based on complex function theory and on eigenprojectors,
as is standard in analytic perturbation theory, is a simplified version of well-known results in the literature.
The second derivation uses a convenient normalization of the right and left eigenvectors
defined in terms of the associated eigenprojector,
but although this dates back to the 1950s, it is rarely discussed in the literature.
We then show how the eigenvector perturbation theory is easily extended to handle other normalizations that are often used in practice.
We also explain how to verify the perturbation results computationally. We conclude with some 
remarks about difficulties introduced by multiple eigenvalues and give references
to work on perturbation of invariant subspaces corresponding to multiple or clustered eigenvalues.
Throughout the paper we give extensive bibliographic commentary and references for further reading.

\end{abstract}

\section{Introduction} \label{intro}
Eigenvalue perturbation theory is an old topic dating originally to the work of Rayleigh in the 19th century.
Broadly speaking, there are two main streams of research.
The most classical is analytic perturbation theory (APT), where one considers the behavior of
eigenvalues of a matrix or linear operator that is an analytic function of one or more parameters.
Authors of well-known books describing this body of work include Kato \cite{Kat66,Kat76, Kat82,Kat95},\footnote{The first edition
of Kato's masterpiece \emph{Perturbation Theory for Linear Operators}
was published in 1966 and a revised second edition appeared in 1976. The most recent edition is the 1995 reprinting
of the second edition with minor corrections. Most of this book is concerned with linear
operators, but the first two chapters treat the finite-dimensional case of matrices, and these appeared as a stand-alone
short version in 1982. Since we are only concerned with matrices in this article, our references to Kato's book are to the 1982
edition, although in any case the equation numbering is consistent across all editions.}
Rellich \cite{Rel69}, Chatelin \cite{Cha11},
Baumg\"artel \cite{Bau85} and, in text book form, Lancaster and Tismenetsky \cite{LanTis85}.
Kato \cite[p.~XII]{Kat82}) and Baumg\"artel \cite[p.\ 21]{Bau85} 
explain that it was Rellich who first established in the 1930s that 
when a Hermitian matrix or self-adjoint linear operator
with an isolated eigenvalue $\lambda$ of multiplicity $m$ is subjected to a real analytic perturbation, 
that is a convergent power series in a real parameter $\kappa$, then
(1) it has exactly $m$ eigenvalues converging to $\lambda$ as $\kappa \to 0$, (2) these eigenvalues can also be expanded in 
convergent power series in $\kappa$ and 
(3) the corresponding eigenvectors can be chosen to be mutually orthogonal and may also be written as convergent power series. 
As Kato notes, these results are exactly what 
were anticipated by Rayleigh, Schr\"odinger and others, but to prove them is by no means trivial, 
even in the finite-dimensional case. 

The second stream of research is largely due to the numerical linear algebra (NLA) community.
It is mostly restricted to matrices and generally concerns perturbation \emph{bounds} rather than expansions, 
describing how to bound the change in the eigenvalues and associated eigenvectors or invariant subspaces when a given
matrix is subjected to a perturbation with a given norm and structure. Here there are a wide
variety of well-known results due to many of the founders of matrix analysis 
and numerical linear algebra: Gerschgorin, Hoffman and Wielandt, Mirksy, Lidskii, Ostrowski, Bauer and Fike, Henrici, Davis and Kahan, Varah, Ruhe, 
Stewart, Elsner, Demmel and others. 
These are discussed in many books, of which the most comprehensive include those by 
Wilkinson \cite{Wil65}, Stewart and Sun \cite{SteSun90}, Bhatia \cite{Bha97,Bha07} and Stewart \cite{Ste01},
as well as Chatelin \cite{Cha12}, which actually covers both the APT and the NLA streams of research in some detail.
See also the survey by Li \cite{Li14}.
An important branch of the NLA stream 
concerns the \emph{pseudospectra} of a matrix;
see the book by Trefethen and Embree \cite{TreEmb05} and the Pseudospectra Gateway web site \cite{EmbTre}.

This paper is inspired by both the APT and the NLA streams of research, and its scope is limited to an important special
case: first-order perturbation analysis of a simple eigenvalue and the corresponding right and left eigenvectors of a 
general square matrix, not assumed to be Hermitian or normal. The eigenvalue result is well known to a broad scientific community. 
The treatment of eigenvectors is more complicated, with a perturbation theory that is not so well known outside
a community of specialists. We give two different proofs of the main eigenvector perturbation theorem. The first, inspired by the NLA research
stream and based on the implicit function theorem, has apparently not appeared in the literature in this form.
The second, based on complex function theory and on eigenprojectors, as is standard in APT,
is largely a simplified version of results in the literature that are well known.
The second derivation uses a convenient normalization of the right and left eigenvectors that depends on the perturbation parameter,
but although this dates back to the 1950s, it is rarely discussed in the literature.
We then show how the eigenvector perturbation theory is easily extended to handle other normalizations that are often used in practice.
We also explain how to verify the perturbation results computationally.
In the final section, we illustrate the difficulties introduced by multiple eigenvalues with two illuminating examples,
and give references to work on perturbation of invariant subspaces corresponding to multiple or clustered eigenvalues.

\section{First-order perturbation theory for a simple eigenvalue} \label{sec:evalpert}

Throughout the paper we use
$\|\cdot\|$ to denote the vector or matrix 2-norm, $I_n$ to denote the identity matrix of order $n$,
the superscript $\T$ to denote transpose and $*$ to denote complex conjugate transpose.
Greek lower case letters denote complex scalars.
Latin lower case letters denote complex vectors, with the exception of $i$ for the imaginary unit and $j,k,\ell,m,n$ for integers.
Upper case letters denote complex matrices or, in some cases, sets in the complex plane.
We begin with an assumption that also serves to establish our notation. 

\begin{assumption}\label{assumption}
Let $A_0 \in \C^{n\times n}$ have a simple eigenvalue $\lambda_0$ corresponding to right eigenvector $x_0$
\emph{(}so $Ax_0=\lambda_0 x_0$ with $x_0\ne 0$\emph{)} and left eigenvector $y_0$ 
\emph{(}so $y_0^*A = \lambda_0 y_0^*$ with $y_0 \ne 0$\emph{)},
normalized so that $y_0^*x_0=1$.
Let $\tau_0\in\C$ and let $A(\tau)$ be a complex-valued matrix function of a complex parameter $\tau$ that is analytic 
in a neighborhood of $\tau_0$, satisfying $A(\tau_0)=A_0$.
\end{assumption}

\begin{remark} \label{remark}
The normalization $y_0^*x_0=1$ is always possible since the right and left eigenvector corresponding to a simple 
eigenvalue cannot be orthogonal. Note that since $x_0$ and $y_0$ are unique only up to scalings, we may multiply 
$x_0$ by any nonzero complex scalar $\omega$ provided we also scale $y_0$ by the reciprocal of the conjugate of 
$\omega$ so that $y_0^*x_0$ remains equal to one. The use of the complex conjugate transpose in $y_0^*$ 
instead of an ordinary  transpose is purely a convention that is often, but not universally, followed. 
The statement that the matrix $A(\tau)$ is
analytic means that each entry of $A(\tau)$ is analytic (equivalently, complex differentiable or holomorphic) in $\tau$ 
in a neighborhood of $\tau_0$.
\end{remark}

The most basic result in eigenvalue perturbation theory follows. 

\begin{theorem} \label{evalpert} {\rm (Eigenvalue Perturbation Theorem)} Under Assumption~\ref{assumption}, 
$A(\tau)$ has a unique eigenvalue $\lambda(\tau)$ that is analytic in a neighborhood of $\tau_0$, 
with $\lambda(\tau_0)=\lambda_0$ and with
\beq\label{lamderiv}
      \lambda'(\tau_0) = y_0^* A'(\tau_0) x_0  
\eeq
where $\lambda'(\tau_0)$ and $A'(\tau_0)$ are respectively the derivatives of $\lambda(\tau)$ 
and $A(\tau)$ at $\tau=\tau_0$.
\end{theorem}
The proof appears in the next section.

The quantity
\[
           \chi = \|x_0\|\|y_0\| \geq |y_0^*x_0| = 1,
\] 
introduced by \cite{Wil65}, is called the \emph{eigenvalue condition number} for $\lambda_0$.
We have $|\lambda'(\tau_0) | \leq \chi \|A'(\tau_0)\|$.
In the real case, $\chi$ is the reciprocal of the cosine of the angle between $x_0$ and $y_0$.
In the special case that $A_0$ is Hermitian, its right and left eigenvectors coincide so $\chi=1$, but in this article
we are concerned with general square matrices.

In the APT research stream, instead of eigenvectors, the focus is mostly on the eigenprojector\footnote{In APT, the standard term is ``eigenprojection", while in NLA, ``spectral projector" is often used. The somewhat nonstandard term 
``eigenprojector" is a compromise.}
corresponding to $\lambda_0$, which can be defined as
\beq \label{Pidef}
          \Pi_0 = x_0 y_0^*
\eeq
and which satisfies
\[
     A_0\Pi_0 = \lambda_0\Pi_0 = \Pi_0A_0 \text{ and }\Pi_0^2 = \Pi_0.
\]  
Note that the eigenprojector does not depend on the normalization used for the eigenvectors 
$x_0$ and $y_0$ (assuming $y_0^*x_0=1$), which simplifies the associated perturbation theory, and note also that $\chi=\|\Pi_0\|$.
Let $\tr$ denote trace and recall the property $\tr(XY) = \tr(YX)$. Clearly, equation \eqref{lamderiv} is equivalent to 
\beq \label{lamderivtrace}     
         \lambda'(\tau_0) = \tr \left( \Pi_0  A'(\tau_0) \right) .
\eeq

Kato \cite[p.\ XIII]{Kat82} explains that the results of Rellich for analytic perturbations of
self-adjoint linear operators were extended (by Sz-Nagy, Kato and others) to non-self-adjoint linear operators and 
therefore non-Hermitian matrices
in the early 1950s using complex function theory, so
\eqref{lamderivtrace}, equivalently \eqref{lamderiv},
was known at that time. However, it seems that these results were not well known
until the publication of the first edition of Kato's book in 1966 
(although Kato did present a summary of these results for the linear case at a conference on matrix computations \cite[p.~104]{Giv58} in 1958).
Eq.\ \eqref{lamderiv} was independently obtained for the analytic case by Lancaster \cite{Lan64}, 
and for the linear case by Wilkinson \cite[p.68--69]{Wil65} and Lidskii \cite{Lid66E}.
They all used the theory of algebraic functions to obtain their results,
exploiting the property that eigenvalues are roots of the characteristic polynomial. A different technique
is used by Stewart and Sun \cite[p.\ 185]{SteSun90} who show that the eigenvalue is differentiable w.r.t.\ its matrix argument using a proof depending on Gerschgorin circles;
the result for a differentiable family $A(\tau)$  then follows from the ordinary chain rule.

We close this section with a brief discussion of multiple eigenvalues. The algebraic multiplicity of $\lambda_0$ is the multiplicity of the factor $\lambda-\lambda_0$ in
the characteristic polynomial $\det(A_0-\lambda I_n)$, while the geometric multiplicity (which is always less than or equal
to the algebraic multiplicity) is the number of associated
linearly independent right (equivalently, left) eigenvectors. A simple eigenvalue has both algebraic and geometric multiplicity equal to one.
More generally, if the algebraic and geometric multiplicity are equal, the eigenvalue is said to be semisimple or nondefective.
An eigenvalue whose geometric multiplicity is one is called nonderogatory.

\section{First-order perturbation theory for an eigenvector corresponding to a simple eigenvalue}  \label{sec:evecpert}

We begin this section with a basic result from linear algebra; see  \cite[Theorem 1.18 and eq.\ (3.10)]{Ste01} for a proof.

\begin{lemma}\label{XYlemma}
Suppose Assumption 1 holds. There exist matrices
$X_1\in \C^{n\times (n-1)}$, $Y_1\in \C^{n\times (n-1)}$ and $B_1 \in \C^{(n-1)\times(n-1)}$ satisfying 
\beq\label{XYeqn}
    X=[x_0,X_1], \quad Y=[y_0,Y_1], \quad  Y^*X = I_n, \quad \text{and} \quad Y^*A_0 X=
                \begin{bmatrix}
                  \lambda_0 & 0 \\
                  0     & B_1
                \end{bmatrix}.
\eeq
\end{lemma}
Note that, from $Y^*X=I_n$, it is immediate that
the columns of $X_1$ and $Y_1$ respectively span the null spaces of $y_0^*$ and $x_0^*$. 
Furthermore, we have
\[
       I_n = XY^* = x_0 y_0^* + X_1 Y_1^*.
\]
We also have $A_0X_1 = X_1 B_1$ and $Y_1^*A_0 = B_1 Y_1^*$, so the columns of $X_1$ and $Y_1$ are respectively bases
for right and left ($n-1$)-dimensional invariant subspaces of $A_0$, and $I_n = \Pi_0 + \Pi_1$ where $\Pi_1=X_1 Y_1^*$
is the complementary projector to $\Pi_0$.
If we assume that $A_0$ is diagonalizable, i.e., with $n$ linearly independent eigenvectors,
then we can take the columns of $X_1$ and of $Y_1$
to respectively be right and left eigenvectors corresponding to the eigenvalues of $A_0$ that differ from
$\lambda_0$, which we may denote by $\lambda_1,\ldots,\lambda_{n-1}$ (some of which could coincide, as
diagonalizability implies only that the eigenvalues are semisimple, not that they are simple). In this case,
we can take $B_1$ to be the diagonal matrix $\diag(\lambda_1,\ldots,\lambda_{n-1})$. More generally, however,
$X_1$ and $Y_1$ may be \emph{any} matrices satisfying \eqref{XYeqn}, 
ignoring the multiplicities and Jordan structure of the other eigenvalues.

Now let
\beq
       S= X_1\big(B_1 - \lambda_0 I_{n-1}\big)^{-1}Y_1^*.   \label{Sdef}
\eeq
It then follows that
\[
         S\Pi_0 = \Pi_0 S = 0 \text{ and } (A_0-\lambda_0 I_n) S = S(A_0-\lambda_0 I_n) = \Pi_1.
\] 
In the NLA stream of research, $S$ is called the \emph{group inverse} of $A_0-\lambda_0 I_n$ \cite[p.\ 240--241]{SteSun90},
\cite{MeySte88}, \cite{CamMey91}, \cite[Theorem 5.2.]{GugOve11}. In the APT research
stream, it is called the \emph{reduced resolvent matrix} of $A_0$ w.r.t.\ the eigenvalue $\lambda_0$
(see \cite[eqs.\ I.5.28 and II.2.11]{Kat82}.)\footnote{The notion of group inverse or reduced resolvent extends beyond
the simple eigenvalue context to multiple eigenvalues. If $\lambda_0$ is a defective eigenvalue with a nontrivial Jordan structure,
the reduced resolvent matrix of $A$ with respect to $\lambda_0$ must take account of ``eigennilpotents".
It is the same as the Drazin inverse of $A_0-\lambda_0 I_n$, a generalization of the group inverse
(see \cite[p.98]{Cha12}, \cite{CamMey91} and, for a method to compute the Drazin inverse, 
\cite{GugOveSte15}).}

We now give a first-order perturbation theorem for right and left eigenvectors corresponding to a simple eigenvalue.

\begin{theorem}\label{evecpert} {\rm (Eigenvector Perturbation Theorem)}
Suppose that Assumption~\ref{assumption} holds and define $X_1$, $Y_1$ and $B_1$ as in Lemma \ref{XYlemma} and $S$ as in \eqref{Sdef}. 
Then there exist  vector-valued functions 
$x( \tau )$ and $y( \tau )^{*}$ that are analytic in a neighborhood of $\tau_0$
with $x(\tau_0)=x_0$, $y(\tau_0)=y_0$ and $y(\tau)^*x(\tau)=1$,
satisfying the right and left eigenvector equations
\beq\label{eveceqns}
A(\tau) x(\tau)=\lambda(\tau) x(\tau), \quad y(\tau)^* A(\tau)=\lambda(\tau) y(\tau)^*
\eeq
where $\lambda(\tau)$ is the analytic function from Theorem~\ref{evalpert}. 
Furthermore, these can be chosen so that 
their derivatives, $x'(\tau)$ and  $(y^*)'(\tau)$, satisfy $y_0^*x'(\tau_0) = 0$ and  $(y^*)'(\tau_0) x_0 = 0$, with\footnote{We use
the notation $(y^*)'(\tau_0)$ to mean $\frac{d}{d \tau} ( y(\tau)^* )|_{\tau=\tau_0}$.}
\begin{align}
x'(\tau_0) & = - S A'(\tau_0) x_0, \label{xderiv} \\
 (y^*)'(\tau_0) & = - y_0^* A'(\tau_0) S. \label{yderiv}
\end{align}

\end{theorem}
Note that it is $y(\tau)^*$, not $y(\tau)$,
that is analytic with respect to the complex parameter $\tau$. However, $y(\tau)$ is
differentiable w.r.t.\ the real and imaginary parts of $\tau$. Note also that we do not claim that $x(\tau)$ and $y(\tau)$ are
unique, even when they are chosen to satisfy \eqref{xderiv} and \eqref{yderiv}. Sometimes, other
normalizations of the eigenvectors, not necessarily satisfying \eqref{xderiv} and \eqref{yderiv},
are preferred, as we shall discuss in \S
\ref{subsec:normalize}. 

It follows from Theorem~\ref{evecpert} that
\[
    \frac{\|x'(\tau_0)\|}{\|x_0\|} \leq  \kappa(X)\|(\lambda_0 I_{n-1}-B_1)^{-1}\| \|A'(\tau_0)\|,
\]
where $\kappa(X)=\kappa(Y)=\|X\|\|Y\|$, the ordinary matrix condition number of $X$, equivalently of $Y$
(as $Y^*=X^{-1}$), with the same bound also holding for $\|(y^*)'(\tau_0)\|/\|y_0\|$.

In the diagonalizable case, as already noted above, we can take $B_1 = \diag(\lambda_1,\ldots,\lambda_{n-1})$, so 
\[
    \frac{ \|x'(\tau_0)\|}{\|x_0\|} \leq  \frac{\kappa(X)\|A'(\tau_0)\|}{\min_{j=1,\ldots,n-1}\{|\lambda_0 - \lambda_j|\}}
 \]
with the same bound also holding for $\|(y^*)'(\tau_0)\|/\|y_0\|$.
In this case, the formula \eqref{xderiv} for the eigenvector derivative was given by
Wilkinson \cite[p.\ 70]{Wil65}. He remarked (p.\  109) that although his derivation is essentially classical 
perturbation theory, a simple but rigorous treatment did not seem to be readily available in the literature.
Lancaster \cite{Lan64} and Lidskii \cite{Lid66E} both showed that the perturbed eigenvector corresponding to a simple eigenvalue 
$\lambda_0$ may be defined to be differentiable at $\lambda_0$, but they did not give the first-order perturbation term. 
The books by Stewart and Sun \cite[sec.\ V.2]{SteSun90} and Stewart \cite[sec.\ 1.3 and 4.2]{Ste01} give 
excellent discussions of the issues summarized above as well as many additional related results.
The eigenvector derivative formula \eqref{xderiv} in Theorem~\ref{evecpert} above is succinctly stated 
just below \cite[eq.\ (3.14), p.\ 46]{Ste01}, where on the same page
Theorem 3.11 stating it more rigorously, and providing additional bounds, is also given;
see also \cite[line 4, p.~48]{Ste01}. The reader is referred to \cite{Ste71} and \cite{Ste73}
for a proof. 
Stewart \cite{Ste71} introduced the idea of establishing the existence of a solution to an algebraic Riccati equation by a
fixed point iteration, a technique that was followed up in \cite[eq.\ (1.5), p.\ 730]{Ste73}
and \cite[eq.\ (7.2), p.187]{Dem86}.  Alternatively, proofs of Theorem~\ref{evecpert} may be derived
by various approaches based on the implicit function theorem; see \cite{Mag85,Sun85} and \cite[Sec.\ 2.1]{Sun98}.
A related argument appears in \cite[Theorem 9.8]{Lax07}. These approaches generally focus on obtaining
results for the right eigenvector subject to some normalization; they can also be applied to obtain results
for the left eigenvector, and these can be normalized further to obtain the condition \mbox{$y(\tau)^*x(\tau)=1$.}
The proof that we give in \S \ref{firstproof} is also based on the implicit function theorem, using a
block-diagonalization approach that obtains the perturbation results for the right and left eigenvectors
simultaneously, ensuring that \mbox{$y(\tau)^*x(\tau)=1$.}  Note, however, that a fundamental difficulty
with eigenvectors is their lack of uniqueness. In contrast, the eigenprojector is uniquely defined, and
satisfies the following perturbation theorem.


\begin{theorem}\label{eprojpert} {\rm (Eigenprojector Perturbation Theorem)}
Suppose Assumption~\ref{assumption} holds and define $X_1$, $Y_1$ and $B_1$ as in Lemma \ref{XYlemma} 
and $S$ as in \eqref{Sdef}. 
Then there exists a matrix-valued function $\Pi ( \tau )$ that is analytic in a neighborhood
of $\tau_0$
with $\Pi(\tau_0)=\Pi_0$, satisfying the eigenprojector equations
\beq \label{eprojeqns}
          A(\tau)\Pi(\tau) = \lambda(\tau)\Pi(\tau) = \Pi(\tau)A(\tau) \text{ and }\Pi(\tau)^2 = \Pi(\tau),
\eeq
and with derivative given by
\beq \label{eigenprojector}
        \Pi'(\tau_0) = -\Pi_0 A'(\tau_0) S - SA'(\tau_0) \Pi_0. 
\eeq
\end{theorem}
This result is well known in the APT research stream \cite[eq.~(II.2.13)]{Kat82}, and, like the eigenvalue perturbation result,
 goes back to the 1950s. Furthermore, while it's easy to see how Theorem~\ref{eprojpert} can be proved using Theorem~\ref{evecpert}, it is also the case that Theorem~\ref{evecpert} can be proved using Theorem~\ref{eprojpert}, by defining the eigenvectors appropriately in terms of the eigenprojector, as discussed below. This provides a convenient
 way to define eigenvectors uniquely.
 
We note that Theorems \ref{evalpert}, \ref{evecpert} and \ref{eprojpert}
 simplify significantly when $A(\tau)$ is a Hermitian function of a real
parameter $\tau$, because then the right and left eigenvectors coincide. The results
for the Hermitian case lead naturally to perturbation theory for singular values and
singular vectors of a general rectangular matrix; see 
\cite[Sec.\ 3.3.1]{Ste01} and \cite[Sec.\ 3.1]{Sun98}.



\subsection{Nonrigorous derivation of the formulas in Theorems~\ref{evalpert}, 
\ref{evecpert}, and \ref{eprojpert}} \label{subsec:nonrigorous}

If we {\em assume} that for $\tau$ in some neighborhood of $\tau_0$, the matrix
$A( \tau )$ has an eigenvalue $\lambda ( \tau )$ and corresponding right and left
eigenvectors $x( \tau )$ and $y ( \tau )$ with
$y( \tau )^{*} x( \tau ) = 1$ such that $\lambda ( \tau )$,
$x( \tau )$ and $y ( \tau )^{*}$ are all analytic functions of $\tau$ satisfying
$\lambda ( \tau_0 ) = \lambda_0$, $x( \tau_0 ) = x_0$, and $y ( \tau_0 ) = y_0$,
then differentiating the equation $A( \tau ) x( \tau ) = \lambda ( \tau ) x( \tau )$
and setting $\tau = \tau_0$, we find
\begin{equation}
A' ( \tau_0 ) x_0 + A_0 x' ( \tau_0 ) = \lambda' ( \tau_0 ) x_0 + \lambda_0 x' ( \tau_0 ) .
\label{nonrigorous1}
\end{equation}
Multiplying on the left by $y_0^{*}$ and using $y_0^{*} A_0 = \lambda_0 y_0^{*}$
and $y_0^{*} x_0 = 1$, we obtain the formula for $\lambda' ( \tau_0 )$:
\[
\lambda' ( \tau_0 ) = y_0^{*} A' ( \tau_0 ) x_0 .
\]

Equation (\ref{nonrigorous1}) can be written in the form
\begin{equation}
( A' ( \tau_0 ) - \lambda' ( \tau_0 ) I_n ) x_0 = - ( A_0 - \lambda_0 I_n ) 
x' ( \tau_0 ) . \label{nonrigorous2}
\end{equation}
Using Lemma \ref{XYlemma}, we can write
\[
A_0 - \lambda_0 I_n = X \left[ \begin{array}{cc} 0 & 0 \\ 0 & B_1 - \lambda_0 I_{n-1}
\end{array} \right] Y^{*} ,
\]
and substituting this into (\ref{nonrigorous2}) and multiplying on the left by $Y^{*}$,
we find
\[
Y^{*} ( A' ( \tau_0 ) - \lambda' ( \tau_0 ) I_n ) x_0 = 
- \left[ \begin{array}{cc} 0 & 0 \\ 0 & B_1 - \lambda_0 I_{n-1}
\end{array} \right] Y^{*} x' ( \tau_0 ) .
\]
The first row equation here is $y_0^{*} ( A' ( \tau_0 ) - \lambda' ( \tau_0 ) I_n ) x_0 = 0$,
which is simply the formula for $\lambda' ( \tau_0 )$.  The remaining $n-1$
equations are
\[
Y_1^{*} ( A' ( \tau_0 ) - \lambda' ( \tau_0 ) I_n ) x_0 = - ( B_1 - \lambda_0 I_{n-1}) Y_1^{*}
x' ( \tau_0 ) ,
\]
and since $Y_1^{*} x_0 = 0$ and $B_1 - \lambda_0 I_{n-1}$ is invertible,
we obtain the following formula for $Y_1^{*} x' ( \tau_0 )$:
\[
Y_1^{*} x' ( \tau_0 ) = - ( B_1 - \lambda_0 I_{n-1} )^{-1} Y_1^{*} A' ( \tau_0 ) x_0 .
\]
Note that $x' ( \tau_0 )$ is not completely determined by this formula because 
each eigenvector $x( \tau )$ is determined only up to a multiplicative constant.
{\em If} we can choose the scale factor in such a way that $y_0^{*} x' ( \tau_0 ) = 0$
then, multiplying on the left by $X_1$ and recalling that $X_1 Y_1^{*} = I_n -
x_0 y_0^{*}$, we obtain the formula in (\ref{xderiv}) for $x' ( \tau_0 )$:
\[
x' ( \tau_0 ) = - X_1 ( B_1 - \lambda_0 I_{n-1} )^{-1} Y_1^{*}  A' ( \tau_0 ) x_0 =
-S A' ( \tau_0 ) x_0 .
\]
Similarly, the formula (\ref{yderiv}) for $( y^{*} )' ( \tau_0 )$ can be derived assuming
that we can choose $y ( \tau )$ so that $( y^{*} )' ( \tau_0 ) x_0 = 0$.

Once formulas (\ref{xderiv}) and (\ref{yderiv}) are established, formula (\ref{eigenprojector})
follows immediately from 
\[
    \left[ x( \tau ) y( \tau )^{*} \right]' = x' ( \tau ) y( \tau )^{*} + x( \tau ) \left[ y( \tau )^{*} \right]' .  
\]
Evaluating at $\tau = \tau_0$ and using formulas (\ref{xderiv}) and (\ref{yderiv}) for $x' ( \tau_0 )$
and $y' ( \tau_0 )$, we obtain formula (\ref{eigenprojector}) for $\Pi' ( \tau_0 )$.

In the following subsections, we establish the {\em assumptions} used here when $\lambda_0$ is
a simple eigenvalue of $A_0$, and thus
obtain {\em proofs} of Theorems \ref{evalpert}, \ref{evecpert}, and \ref{eprojpert},
in two different ways.
The first involves finding equations that a similarity transformation must satisfy if 
it is to take $A( \tau )$ (or, more specifically, $Y^{*} A( \tau ) X$) to a block diagonal form like that 
in Lemma \ref{XYlemma} for $A_0$.  The
implicit function theorem\footnote{Since the perturbation parameter $\tau$ and the matrix family $A(\tau)$ are complex, we need
a version of the implicit function theorem from complex analysis, but in the special case that $\tau$ and $A(\tau)$ are real, we could use a more familiar version
from real analysis. In that case, although some of the eigenvalues and eigenvectors of a real matrix may be complex, they occur in complex conjugate pairs
and are easily represented using real quantities.} is then invoked to show that these equations have a unique solution, for $\tau$ in
some neighborhood of $\tau_0$, and that the solution is analytic in $\tau$.  The second uses 
the argument principle and the residue theorem from complex analysis to establish that, for $\tau$ in a neighborhood of $\tau_0$,
each matrix $A( \tau )$ has a simple eigenvalue $\lambda ( \tau )$ that is analytic in $\tau$ and 
satisfies $\lambda ( \tau_0 ) = \lambda_0$.  It then follows from Lemma~\ref{XYlemma} that there is
a similarity transformation taking $A( \tau )$ to block diagonal form, but Lemma~\ref{XYlemma} says
nothing about analyticity or even continuity of the associated matrices $X( \tau )$ and $Y( \tau )^{*}$.
Instead, the similarity transformation is applied to the {\em resolvent} and integrated to 
obtain an expression for the eigenprojector $\Pi ( \tau )$ that is shown to be analytic in $\tau$.
Finally, left and right eigenvectors satisfying the analyticity conditions along
with the derivative formulas (\ref{xderiv}) and (\ref{yderiv}) are defined in terms of the eigenprojector.

Note that the assumptions used here do \emph{not} generally hold when $\lambda_0$ is not a simple eigenvalue of $A_0$,
as discussed in \S \ref{multiple}.

\subsection{First proof of Theorems~\ref{evalpert}, \ref{evecpert} and \ref{eprojpert}, using techniques from the NLA research stream}
\label{firstproof}

The first proof that we give is inspired by the NLA research steam, but 
instead of Stewart's fixed-point iteration technique mentioned previously, we rely
on the implicit function theorem \cite[Theorem 1.4.11]{Kra01}, which we now state in the form that we need.

\begin{theorem}\label{holift} {\rm (Implicit Function Theorem)}
Let $D \subset \mathbb{C} \times \mathbb{C}^{\ell}$ be an open set, $h = ( h_1 , \ldots ,
h_{\ell} ) : D \rightarrow \mathbb{C}^{\ell}$ an analytic
 mapping, and $( \tau_0 , z^{0} ) \in D$
a point where $h ( \tau_0 , z^{0} ) = 0$ and where the Jacobian matrix $\left( \frac{\partial h_j}
{\partial z_k} \right)_{j,k=1}^{\ell}$ is nonsingular.  Then the system of equations $h( \tau, z ) = 0$
has a unique analytic solution $z = z( \tau )$ in a neighborhood of $\tau_0$ that satisfies $z( \tau_0 ) = z^0$.
\end{theorem}


We now exploit this result in our proof of Theorem \ref{evecpert}.
The setting of the stage before applying the implicit function theorem follows Demmel's variant of 
Stewart's derivation mentioned above.
We obtain a proof of Theorem \ref{evalpert} along the way, and then give a proof of Theorem \ref{eprojpert} as an easy consequence.

Using Lemma \ref{XYlemma}, define
\beq \label{YAXeqn}
\begin{bmatrix}
             \gamma_{11}(\tau) & c_{12}^*(\tau)  \\
               c_{21}(\tau) & C_{22}(\tau)
    \end{bmatrix}  =
Y^{*}A(\tau)X   
   =\begin{bmatrix}
        \lambda_0 & 0 \\
                  0     & B_1
    \end{bmatrix} + Y^{*}(A(\tau) - A_0)X.
\eeq
Here the scalar $\gamma_{11}$, the row and column vectors $c_{12}^*$ and $c_{21}$ and
the $(n-1)\times (n-1)$ matrix $C_{22}$ are analytic functions of $\tau$ near $\tau_0$, since $A( \tau )$ is.
In what follows, we will transform this matrix into a block diagonal matrix by a similarity transformation.
We will choose $p(\tau)$, $q(\tau)$ so that
\[
P(\tau)=\begin{bmatrix}
    1 & -q(\tau) ^{*}\\
    p(\tau) & I_{n-1}
  \end{bmatrix}, \quad
Q(\tau)=\begin{bmatrix}
    1 & -p(\tau)^{*} \\
    q(\tau) & I_{n-1}
  \end{bmatrix},
\]
\[
D(\tau)=Q(\tau)^{*}P(\tau)=\begin{bmatrix}
          1+q(\tau)^{*}p(\tau) & 0 \\
          0 & I_{n-1}+p(\tau)q(\tau)^{*}
        \end{bmatrix},
\]
with $p(\tau)$ and $q(\tau)^{*}$, and consequently $P(\tau)$, $Q(\tau)^{*}$, and $D(\tau)$, analytic
in a neighborhood of $\tau_0$, with $p(\tau_0) = q(\tau_0) = 0$, and hence $P(\tau_0)=Q(\tau_0)=D(\tau_0)=I_n$.
This transformation idea traces back to \cite[p.730]{Ste73} who designed a $P$ with $p=q$, but
in the form given here, it is due to \cite[p.187]{Dem86}.

We would like to have, for $\tau$ sufficiently close to $\tau_0$, the similarity transformation
\beq \label{blkdiag}
Q(\tau)^{*}Y^{*}A(\tau)XP(\tau)D(\tau)^{-1}=\begin{bmatrix}
                       \lambda(\tau) & 0 \\
                       0 & B(\tau)
                     \end{bmatrix},
\eeq
where $\lambda(\tau)$ and $B(\tau)$ are also analytic, with $\lambda(\tau_0)=\lambda_0$ and
$B(\tau_0)=B_1$. Since $D(\tau)$ is block diagonal by definition, we need $Q(\tau)^{*}Y^{*}A(\tau)XP(\tau)$ to be block
diagonal.  Suppressing the dependence on $\tau$, 
this last matrix is given by
\beq \label{blkdiag2}
\begin{bmatrix}
     \gamma_{11}+q^{*}c_{21}+c_{12}^*p+q^{*}C_{22}p & -q^{*}(\gamma_{11}I_{n-1}-C_{22})+c_{12}^*-q^{*}c_{21}q^{*} \\
     -(\gamma_{11}I_{n-1}-C_{22})p+c_{21}-pc_{12}^*p & C_{22}-pc_{12}^*-c_{21}q^{*}+p\gamma_{11}q^{*}
   \end{bmatrix}.
\eeq
For clarity, we introduce the notation $w(\tau)$ for the analytic row vector function $q(\tau)^{*}$. We then
seek column and row vector analytic functions $p(\tau)$ and $w(\tau)$ making the off-diagonal blocks of \eqref{blkdiag2}
zero, i.e., satisfying 
\begin{align}
f(\tau,p(\tau))&:=-\left [\gamma_{11}(\tau)I_{n-1}-C_{22}(\tau)\right ]p(\tau)+c_{21}(\tau)-p(\tau)c_{12}(\tau)^*p(\tau)=0, \label{feqn}\\
g(\tau,w(\tau))&:=-w(\tau)\left [\gamma_{11}(\tau)I_{n-1}-C_{22}(\tau) \right ] +c_{12}(\tau)^*-w(\tau)c_{21}(\tau)w(\tau)=0. \label{geqn}
\end{align}

Taking $\ell = n-1$, $z = 0$, and $h$ equal to first $f$ and then $g$ with $s$
equal to $p$ and $w$, respectively, in Theorem \ref{holift}, we note that since
$f(\tau_0,0)=0$, $g(\tau_0,0)=0$, and the Jacobian matrices
$$
\left.\frac {\partial f(\tau,p)}{\partial p}\right|_{(\tau,p)=(\tau_0,0)}=-(\lambda_0 I_{n-1}-B_1), \quad
\left.\frac {\partial g(\tau,w)}{\partial w}\right|_{(\tau,w)=(\tau_0,0)}=-(\lambda_0 I_{n-1}-B_1),
$$
are nonsingular, 
there are unique functions $p$ and $w$, analytic in a neighborhood of $\tau_0$, 
satisfying \eqref{feqn} and \eqref{geqn} with $p(\tau_0)=0$, $w(\tau_0)=0$, and
\beq \label{pwderivs}
      p'(\tau_0)=(\lambda_0 I_{n-1} -B_1)^{-1}c_{21}'(\tau_0), \quad w'(\tau_0)=(c_{12}^*)'(\tau_0)(\lambda_0 I_{n-1}-B_1)^{-1},
\eeq
where, using the definition \eqref{YAXeqn}, we have
\beq \label{a21a12deriv}
      c_{21}'(\tau_0) = Y_1^*A'(\tau_0) x_0 \quad \text{and} \quad (c_{12}^*)'(\tau_0) = y_0^*A'(\tau_0) X_1.
\eeq
Thus, \eqref{blkdiag2} is block diagonal and hence \eqref{blkdiag} holds, with $\lambda(\tau)$ an eigenvalue of $A(\tau)$,
and with
\begin{align*}
\lambda(\tau)&=[\gamma_{11}+wc_{21}+c_{12}^*p+wC_{22}p](1+wp)^{-1}, \\
B(\tau)&=[C_{22}-pc_{12}^*-c_{21}w+p\gamma_{11}w](I_{n-1}+pw)^{-1},
\end{align*}
again suppressing dependence on $\tau$ on the right-hand sides. These functions are analytic in a neighborhood of
$\tau_0$, satisfying $\lambda(\tau_0)=\lambda_0$ and $B(\tau_0)=B_1$,
with 
\[
    \lambda'(\tau_0) = \gamma_{11}'(\tau_0) = y_0^*A'(\tau_0)x_0,
\]
proving Theorem~\ref{evalpert}. 

Let $e_1$ be the first column of the identity matrix $I_n$. Multiplying \eqref{blkdiag} on the left by $X Q(\tau)^{-*}$
and on the right by $e_1$ we obtain, using $Y^*X=I$ and $D(\tau)=Q(\tau)^{*}P(\tau)$,
\[ 
    A(\tau)X P(\tau)D(\tau)^{-1} e_1 = \lambda(\tau) X P(\tau)D(\tau)^{-1}  e_1,
\]  
so
\[
       x(\tau) :=  X P(\tau)D(\tau)^{-1} e_1 = \frac{1}{1+q(\tau)^* p(\tau)}\left [ x_0 + X_1 p(\tau) \right ]
\]
is analytic with $x(\tau_0)=x_0$ and satisfies the desired right eigenvector equation in \eqref{eveceqns}.
Likewise, multiplying \eqref{blkdiag} on the
right by $D(\tau)P(\tau)^{-1}Y^*$ and on the left by $e_1^{\T}$ gives
\[
        e_1^{\T}Q(\tau)^*Y^* A(\tau) = \lambda(\tau) e_1^{\T} Q(\tau)^* Y^*
\]
so
\[      
     y(\tau)^* := e_1^{\T}Q(\tau)^*Y^* = y_0^* + q(\tau)^* Y_1^*
\]
is analytic with $y(\tau_0)=y_0$ and satisfies the left eigenvector equation in \eqref{eveceqns}.
Furthermore, $y(\tau)^*x(\tau)=1$, as claimed.
Finally, differentiating $x(\tau)$ and $y(\tau)^*$ we have
\[
        x'(\tau_0) = X_1p'(\tau_0) \quad \text{and} \quad (y^*)'(\tau_0) = (q^*)'(\tau_0) Y_1^*
\]
so combining these with \eqref{pwderivs} and
\eqref{a21a12deriv}, recalling that $w(\tau)=q(\tau)^*$, and using the definition of $S$ in \eqref{Sdef}, we obtain the eigenvector
derivative formulas \eqref{xderiv}, \eqref{yderiv}. The properties $y_0^*x'(\tau_0) = 0$ and  $(y^*)'(\tau_0) x_0 = 0$ follow,
so Theorem~\ref{evecpert} is proved. 

Finally, define
\[
          \Pi(\tau)=x(\tau) y(\tau)^*.
\]
The eigenprojector equations \eqref{eprojeqns} follow immediately. We have
\[
           \Pi'(\tau_0) = x'(\tau_0) y(\tau_0)^* + x(\tau_0) (y^*)'(\tau_0),
\]
so Theorem~\ref{eprojpert} follows from \eqref{xderiv} and \eqref{yderiv}.

\hfill$\Box$

\subsection{Second proof of Theorems~\ref{evalpert}, \ref{evecpert} and \ref{eprojpert}, using techniques from the APT research stream}
\label{secondproof}

In this proof, in contrast to the previous one, we focus on proving Thereom~\ref{eprojpert} first, 
obtaining the proof of Theorem~\ref{evalpert} along
the way, and finally obtaining Theorem~\ref{evecpert} as a consequence.
This proof of Theorem~\ref{eprojpert} is based on complex function theory, as is standard in APT. 
However, our derivation is simpler than most given in the literature, which usually prove more
general results, such as giving complete analytic expansions for the eigenvalue and eigenprojector, while we are concerned
only with the first order term. The key to the last part of the proof, yielding
Theorem~\ref{evecpert}, is to use an appropriate eigenvector normalization.

The main tool here is the residue theorem \cite[p.~293-294, Thm.~8.1 and 8.2]{MatHow06}:
\begin{theorem}\label{residuethm} {\rm (Residue Theorem)} 
Let $D$ be a simply connected domain in $\C$ and let $\Gamma$ be a simple closed positively oriented
contour that lies in $D$.  If $f$ is analytic inside $\Gamma$ and on $\Gamma$, except
at the points $\zeta_1 , \ldots , \zeta_m$ that lie inside $\Gamma$, then
\[
\int_{\Gamma} f(\zeta)\,d\zeta = 2 \pi i \sum_{\ell =1}^m \Res [ f, \zeta_{\ell} ] ,
\]
where if $f$ has a simple pole at $\zeta_{\ell}$, then
\[
\Res [ f, \zeta_{\ell} ] = \lim_{\zeta \rightarrow \zeta_{\ell}} (\zeta - \zeta_{\ell} ) f(\zeta) ,
\]
and if $f$ has a pole of order $k$ at $\zeta_{\ell}$, then
\[
\Res [ f, \zeta_{\ell} ] = \frac{1}{(k-1)!} \lim_{\zeta \rightarrow \zeta_{\ell}} \frac{d^{k-1}}{d \zeta^{k-1}} \left [ (\zeta - \zeta_{\ell} )^k f(\zeta) \right ] .
\]
\end{theorem}

Let $\Gamma$ be the boundary of an open set
$\Delta$ in the complex plane containing the simple eigenvalue 
$\lambda_0$, with no other eigenvalues of $A_0$ in $\Delta \cup \Gamma$.
First note that since $\det(A_0-\lambda I_n)$, the characteristic polynomial $p_0$ of $A_0$, 
does not vanish on $\Gamma$, the same will hold for all polynomials with coefficients
sufficiently close to those of $p_0$; in particular, it will hold for $\det(A(\tau)-\lambda I_n)$, the characteristic polynomial
$p_{\tau}$ of $A( \tau )$, if $\tau$ is sufficiently close to $\tau_0$; say, $| \tau - \tau_0 | \leq \epsilon$.
From here on, we always assume that $| \tau - \tau_0 | \leq \epsilon$.
By the argument principle \cite[p.~328, Thm.~8.8]{MatHow06},
the number of zeros of $p_{\tau}$ inside $\Delta$ is
\[
\frac{1}{2 \pi i} \int_{\Gamma} \frac{\frac{d}{d\zeta}p_{\tau} (\zeta)}{p_{\tau} (\zeta)}\,d\zeta .
\]
For $\tau = \tau_0$, this value is $1$.  Since
for each $\zeta \in \Gamma$,  the integrand $\frac{d}{d\zeta}p_{\tau} (\zeta)/p_{\tau} (\zeta)$ is a continuous function of $\tau$,
the integral above is as well.  Since it is integer-valued, it must be the constant $1$.
So, let $\lambda(\tau)$ denote the unique root of
$p_{\tau}$ in the region $\Delta$, i.e., the unique eigenvalue of $A( \tau )$ in $\Delta$.
Note that this means that $p_{\tau} (\zeta)$ can be written in the form $(\zeta - \lambda ( \tau )) q(\zeta)$, where
$q$ has no roots in $\Delta \cup \Gamma$.  It therefore follows from the residue theorem that
\begin{equation}
\frac{1}{2 \pi i} \int_{\Gamma} \frac{\frac{d}{d\zeta}p_{\tau} (\zeta)}{p_{\tau} (\zeta)} \zeta\,d\zeta = \lambda( \tau )
\label{lambdaoftau}
\end{equation}
since
\[
\Res \left[ \frac{\frac{d}{d\zeta} ( p_{\tau} (\zeta) )~\zeta}{p_{\tau} (\zeta)} , \lambda ( \tau ) \right] =
\lim_{\zeta \rightarrow \lambda ( \tau )} \frac{\left( \zeta - \lambda ( \tau ) \right) 
\left[ (\zeta - \lambda ( \tau )) \frac{d}{d\zeta} q(\zeta) + q(\zeta) \right] \zeta} {\left( (\zeta - \lambda ( \tau ) \right) q(\zeta)}
= \lambda ( \tau ) .
\]
Since the left-hand side of (\ref{lambdaoftau}) is an analytic function of $\tau$, the right-hand side is as well.
Thus $A( \tau )$ has a unique eigenvalue $\lambda ( \tau )$ in $\Delta$ and $\lambda ( \tau )$ is an
analytic function of $\tau$.

For $\zeta$ not an eigenvalue of $A ( \tau )$, define
the {\em resolvent} of $A( \tau )$ by
\begin{equation}
        R( \zeta ; A( \tau ) ) := ( A( \tau ) - \zeta I_n )^{-1} . \label{resolvent}
\end{equation}
Lemma \ref{XYlemma} states that there exist left and right eigenvectors 
$y_0 ( \tau )$ and $x_0 ( \tau )$ associated with $\lambda ( \tau )$ and satisfying
$y_0 ( \tau )^* x_0 ( \tau ) = 1$, along with matrices $X_1 ( \tau ) \in 
\mathbb{C}^{n \times (n-1)}$, $Y_1 ( \tau ) \in \mathbb{C}^{n \times (n-1)}$ and
$B_1 ( \tau ) \in \mathbb{C}^{(n-1) \times (n-1)}$, satisfying
\[
X( \tau ) = [ x_0 ( \tau ) , X_1 ( \tau ) ] ,~~
Y( \tau ) = [ y_0 ( \tau ) , Y_1 ( \tau ) ] ,~~
Y( \tau )^{*} X( \tau ) = I_n ,~~\mbox{and}~~
\]
\[
Y( \tau )^{*} A( \tau ) X ( \tau ) = \left[ \begin{array}{cc} \lambda ( \tau ) & 0 \\
0 & B_1 ( \tau ) \end{array} \right] .
\]
Note that we do not claim that $X(\tau)$ and $Y(\tau)^*$ are analytic, or even continuous functions of $\tau$. 
It follows that the resolvent of $A( \tau )$ satisfies
\begin{eqnarray}
R( \zeta; A( \tau ) ) & = & X ( \tau ) \left[ \begin{array}{cc} 
( \lambda ( \tau ) - \zeta )^{-1} & 0 \\ 
0 & ( B_1 ( \tau ) - \zeta I_{n-1} )^{-1} \end{array} \right] Y( \tau )^{*} \nonumber \\
 & = & ( \lambda ( \tau ) - \zeta )^{-1} x_0 ( \tau ) y_0 ( \tau )^{*} + 
S( \zeta ; \tau ) , \label{partialfractions}
\end{eqnarray}
where $S( \zeta; \tau ) = X_1 ( \tau ) ( B_1 ( \tau ) - \zeta I_{n-1} )^{-1} 
Y_1 ( \tau )^{*}$.
Now, $S( \zeta ; \tau )$ is a matrix-valued function of $\zeta$ with no poles in $\Delta$,
so it follows from the residue theorem (applied to the functions associated with each entry of $S(\zeta; \tau )$) that 
$\int_{\Gamma} S( \zeta; \tau )\,d\zeta = 0$ and therefore from (\ref{partialfractions}) that
\begin{eqnarray}
- \frac{1}{2 \pi i} \int_{\Gamma} R( \zeta ; A( \tau ) )\,d\zeta & = &  
- \frac{1}{2 \pi i} \left( \int_{\Gamma} ( \lambda ( \tau ) - \zeta )^{-1}\,d\zeta \right) x_0 ( \tau ) y_0 ( \tau )^{*} 
\nonumber \\
& = & x_0 ( \tau ) y_0 ( \tau )^{*} = \Pi ( \tau ) . \label{Pi_tau}
\end{eqnarray}
For $\tau = \tau_0$, this is $\Pi ( \tau_0 ) = \Pi_0$.

From the definition of the resolvent (\ref{resolvent}), it follows that since
$A( \tau )$ is an analytic function of $\tau$, the resolvent $R(\zeta; A( \tau ))$
is as well, provided that $\zeta$ is not an eigenvalue of $A( \tau )$.
Differentiating the equation 
\[
     (A( \tau ) - \zeta I_n ) R( \zeta; A( \tau ) ) = I_n
\]
with respect to $\tau$ gives
\[
\frac{\partial}{\partial \tau} R( \zeta; A( \tau ) ) = - R( \zeta; A( \tau ) ) A' ( \tau )
R( \zeta; A( \tau ) ) .
\]
Considering expression (\ref{Pi_tau}) for $\Pi(\tau)$, it follows that $\Pi(\tau)$
is also analytic, and its derivative at $\tau = \tau_0$ is
\[
\Pi' ( \tau_0 ) =  \frac{1}{2 \pi i} \int_{\Gamma} R( \zeta ; A_0 ) A' ( \tau_0 )
R( \zeta ; A_0 )\,d\zeta .
\]
Using (\ref{partialfractions}), and writing $S_0( \zeta )$ for $S( \zeta ; \tau_0 )$, we obtain
\begin{eqnarray*}
\Pi' ( \tau_0 ) & = & \frac{1}{2 \pi i} \int_{\Gamma} \left[ ( \lambda_0 -\zeta)^{-1} \Pi_0 + 
S_0 ( \zeta ) \right] A'( \tau_0 ) \left[ ( \lambda_0 - \zeta )^{-1} \Pi_0 + S_0 ( \zeta ) \right]\,d\zeta \\
 & = & \frac{1}{2 \pi i} \left( \int_{\Gamma} ( \lambda_0 - \zeta )^{-2}\,d\zeta \right) \Pi_0 A' ( \tau_0 ) \Pi_0
+ \frac{1}{2 \pi i} \int_{\Gamma} S_0(\zeta ) A' ( \tau_0 ) S_0(\zeta)\,d\zeta  \\
 &  & + \ \frac{1}{2 \pi i} \int_{\Gamma}  ( \lambda_0 - \zeta)^{-1} \left [ \Pi_0 A' ( \tau_0 )
S_0 ( \zeta ) + S_0 ( \zeta ) A' ( \tau_0 ) \Pi_0  \right]\,d\zeta .  
\end{eqnarray*}
From the residue theorem, the first term is zero since the integrand has a pole of order 2 at $\lambda_0$,
with $\Res [( \lambda_0 - \zeta )^{-2},\lambda_0] = 0$.
The second term is also zero because the
integrand has no poles inside $\Gamma$.
Since the integrand of the remaining term has a simple pole at $\lambda_0$ with
$\Res[(\lambda_0-\zeta)^{-1},\lambda_0]=-1$, we have
\begin{equation}
\Pi' ( \tau_0 ) = - \Pi_0 A' ( \tau_0 ) S_0( \lambda_0 ) - S_0( \lambda_0 ) A' ( \tau_0 ) \Pi_0 ,
\label{Piprime}
\end{equation}
where $S_0( \lambda_0 ) = S(\lambda_0;\tau_0) = X_1 ( \tau_0 ) ( B_1 ( \tau_0 ) - \lambda_0 I )^{-1} 
Y_1 ( \tau_0 )$ is the same as $S$ defined in (\ref{Sdef}).
This proves Theorem~\ref{eprojpert}.

Now we define the eigenvectors in terms of the eigenprojector, using a normalization that goes back to
\cite{SzNa51} (see \cite[eq.\ (7.1.12)]{Bau85}): 
\beq \label{xytaudef}
x( \tau ) =  (y_0^* \Pi(\tau)x_0)^{-1/2} \Pi(\tau) x_0, \quad
y( \tau )^* = (y_0 ^*\Pi(\tau)x_0)^{-1/2}  y_0^* \Pi(\tau) ,
\eeq
where we use the principal branch of the square root function (and assume that $\epsilon$ is small enough
so that $| \tau - \tau_0 | \leq \epsilon$ implies that the quantities under the square roots, which are $1$
for $\tau = \tau_0$, are bounded away from the origin and the negative real axis).\footnote{See also
\cite[eq.\ (II.3.24)]{Kat82}, which uses a related definition but without the square root, 
resulting in $y(\tau)^*x_0=y_0^*x(\tau)=1$ instead of $y(\tau)^*x(\tau) = 1$.}
Since $\Pi ( \tau )$ is analytic, it follows that $x( \tau )$ and $y( \tau )^*$ are as well. From (\ref{xytaudef}) we have
\[
       A(\tau) x(\tau) =  (y_0^* \Pi(\tau)x_0)^{-1/2} A(\tau)  \Pi(\tau) x_0 =  (y_0^* \Pi(\tau)x_0)^{-1/2} \lambda(\tau)\Pi(\tau) x_0 = \lambda(\tau) x(\tau)
\]
and similarly $y(\tau)^*A(\tau) = \lambda(\tau) y(\tau)^*$, so the eigenvector equations \eqref{eveceqns} hold as required,
and, since $\Pi(\tau)^2 = \Pi(\tau)$, we obtain
\[
    y(\tau)^*x(\tau) = (y_0^* \Pi(\tau)x_0)^{-1} y_0^*\Pi(\tau)^2 x_0 = 1
\]
as claimed in Theorem~\ref{evecpert}.

To obtain the eigenvalue derivative, we differentiate the equation
\[
      (A( \tau ) - \lambda ( \tau ) I_n ) x( \tau ) = 0
\]
and evaluate it at $\tau = \tau_0$:
\[
( A' ( \tau_0 ) - \lambda' ( \tau_0 ) I_n ) x_0 + ( A_0  - \lambda_0 I_n ) x' ( \tau_0 ) = 0 .
\]
Multiplying by $y_0^{*}$ on the left, this becomes
\[
\lambda' ( \tau_0 ) = y_0^{*} A' ( \tau_0 ) x_0,
\]
proving Theorem \ref{evalpert}.

Finally, using \eqref{eigenprojector}, we have
\begin{align*}
      x'(\tau_0) & = -\frac{1}{2} (y_0^* \Pi_0 x_0)^{-3/2} (y_0^* \Pi'(\tau_0) x_0) \Pi_0 x_0  +   [y_0^* \Pi_0 x_0]^{-1/2} \Pi'(\tau_0) x_0\\
       & = -\frac{1}{2} ( y_0^*[-\Pi_0 A'(\tau_0) S - S A'(\tau_0) \Pi_0] x_0 ) \Pi_0 x_0  - [\Pi_0 A'(\tau_0) S +S A'(\tau_0) \Pi_0]  x_0.
\end{align*}
The first three terms are zero since $Sx_0 = -X_1(\lambda_0 I - B_1)^{-1} Y_1^*x_0 = 0$ and likewise $y_0^* S = 0$. So,
as $\Pi_0x_0 = x_0$, we have
\[
      x'(\tau_0) =  -S A'(\tau_0)x_0.
\]
Similarly, $(y^*)'(\tau_0) =- y_0^* A'(\tau_0) S$. The properties $y_0^*x'(\tau_0) = 0$ and  
$(y^*)'(\tau_0) x_0 = 0$ follow,
so the proof of Theorem~\ref{evecpert} is complete.

\hfill$\Box$

\subsection{Eigenvector normalizations}
\label{subsec:normalize}

Theorem~\ref{evecpert} does not state formulas for the analytic eigenvector functions $x(\tau)$ and $y(\tau)^*$, specifying only that they exist, satisfying $y(\tau)^*x(\tau)=1$, with derivatives
given by \eqref{xderiv} and \eqref{yderiv}. Furthermore, the first proof given in \S \ref{firstproof} does not provide formulas
for $x(\tau)$, $y(\tau)$, showing only that they exist via the implicit function theorem. However, the second proof given in \S \ref{secondproof}
\emph{does} provide formulas for $x(\tau)$ and $y(\tau)$ in terms of the eigenprojector $\Pi(\tau)$ (which \emph{is} uniquely defined)
and the eigenvectors $x_0$ and $y_0$.  This formula \eqref{xytaudef} may be viewed as a \emph{normalization}  because it
provides a way to define the eigenvectors uniquely, and furthermore, it has the property that $x(\tau)$ and $y(\tau)^*$ are analytic near $\tau_0$ and 
satisfy $y(\tau)^*x(\tau) = 1$. Let us refer to \eqref{xytaudef} as ``normalization 0".

Although the beautifully simple normalization \eqref{xytaudef} dates to the 1950s, it seems to be rarely used. 
In this subsection we discuss some other normalizations that are more commonly
used in practice. Let us denote the resulting normalized eigenvectors by $\hat x(\tau)$ and $\hat y(\tau)$,
and relate them to $x(\tau)$ and $y(\tau)$, as defined in \eqref{xytaudef}, by
\begin{equation}\label{xhatyhat}
\hat x(\tau)=\alpha(\tau)\,x(\tau), \quad  \hat y(\tau)^*=\beta(\tau)^*\,y(\tau)^*,  
\end{equation}
where $\alpha(\tau)$ and $\beta(\tau)$ are two nonzero complex-valued scalar functions of $\tau$  to be defined below.
Here we use the complex conjugate
of $\beta$ in the definition to be consistent with the conjugated left eigenvector notation.
The analyticity of $\hat x(\tau)$ and $\hat y(\tau)^{*}$  near $\tau_0$ depends on that of $\alpha(\tau)$ and $\beta(\tau)^*$.
We consider several possible normalizations, continuing to assume 
that $y_0^*x_0 = 1$ but not necessarily that $\hat{y} ( \tau )^{*} \hat{x} ( \tau )= 1$ for $\tau\not = 0$. In all cases, the formula 
\beq \label{xhatyhat_derivs}
    \hat x'(\tau_0) =  \alpha(\tau_0)x'(\tau_0) + \alpha'(\tau_0)x_0 , \quad 
    (\hat y^*)'(\tau_0) =  \beta^*(\tau_0)(y^*)'(\tau_0) + (\beta^*)'(\tau_0)y_0^* 
\eeq
follows immediately from \eqref{xhatyhat}, so to determine the derivatives $\hat x'(\tau_0)$ and $(\hat y^*)'(\tau_0)$,
we need only determine $\alpha(\tau_0)$,
$\alpha'(\tau_0)$, $\beta^*(\tau_0)$ and $(\beta^*)'(\tau_0)$,
obtaining $x'(\tau_0)$ and $(y^*)'(\tau_0)$ from \eqref{xderiv}, \eqref{yderiv}, as stated in Theorem~\ref{evecpert}.
Note that the derivatives of the normalized eigenvectors, $\hat x'(\tau_0)$ and $ (\hat y^*)'(\tau_0)$, 
do \emph{not} necessarily satisfy $y_0^*\hat x'(\tau_0) = 0$ and  $(\hat y^*)'(\tau_0) x_0 = 0$,
unlike the derivatives $x'(\tau_0)$ and $(y^*)'(\tau_0)$.
We now define several different normalizations:

\begin{enumerate}
  \item $e_1^{\T}\hat x(\tau)=e_1^{\T}\hat y(\tau)=1$ (i.e., the first entries of $\hat x(\tau)$ and $\hat y(\tau)$ are one). 
  This is possible for $\tau$
        sufficiently close to $\tau_0$ if $e_1^{\T} x_0\ne 0$ and $e_1^{\T} y_0\ne 0$. 
        Suppose this is the case. Then 
\beq \label{alphatau1}
      \alpha(\tau) = \frac{1}{e_1^{\T} x(\tau)} \quad \text{so}\quad \alpha'(\tau_0) = -\frac{e_1^{\T}x'(\tau_0)}{(e_1^{\T}x_0)^2} 
\eeq
and 
\beq \label{betatau1}
      \beta(\tau)^* = \frac{1}{y(\tau)^*e_1} \quad \text{so}\quad 
     (\beta^*)'(\tau_0) = -\frac{(y^*)'(\tau_0)e_1}{(y_0^* e_1)^2} .
\eeq
Here $\alpha(\tau)$, $\beta(\tau)^*$ and hence $\hat x(\tau)$, $\hat y^*(\tau)$ are analytic near $\tau_0$.

 \item $e_1^{\T}\hat x(\tau)=1$ and $\hat y(\tau)^{*}\hat x(\tau)=1$. This normalization is defined near $\tau_0$
 under the assumption that $e_1^{\T} x_0\ne 0$; no additional assumption is needed since $y_0^{*}x_0 = 1$.
Clearly \eqref{alphatau1} holds as before. In addition, we have
 \beq \label{betatau2}
        \beta(\tau)^*=\frac{e_1^{\T}x(\tau)}{y(\tau)^*x(\tau)} = e_1^{\T}x(\tau), \quad
         (\beta^*)'(\tau_0) = e_1^{\T} x'(\tau_0) .
 \eeq
Again, $\alpha(\tau)$ and $\beta(\tau)^*$ are analytic near $\tau_0$. 
        
\item $\hat x(\tau)^{\T}\hat x(\tau)=1$ and $\hat y(\tau)^{*}\hat x(\tau)=1$. This normalization is defined near $\tau_0$
if $x_0^{\T}x_0\ne 0$, which may not be the case when $x_0$ is complex. Suppose this does hold. We have
\beq \label{alphatau3}
   \alpha(\tau) = \pm \frac{1}{(x(\tau)^{\T}x(\tau))^{1/2}}, \quad \alpha'(\tau_0)=\mp \frac{x'(\tau_0)^{\T}x_0}{(x_0^{\T}x_0)^{3/2}},
\eeq
\beq  \label{betatau3}
      \beta(\tau)^* = \pm\frac{(x(\tau)^{\T}x(\tau))^{1/2}}{y(\tau)^*x(\tau)} = \pm(x(\tau)^{\T}x(\tau))^{1/2}, 
       (\beta^*)'(\tau_0)=\pm \frac{x'(\tau_0)^{\T}x_0}{(x_0^{\T}x_0)^{1/2}}. 
\eeq
Either sign may be used, resulting in analytic $\alpha(\tau)$ and $\beta^*(\tau)$ near $\tau=\tau_0$. 
           
\item $\hat x(\tau)^*\hat x(\tau) = 1$ and $\hat y(\tau)^{*}\hat x(\tau)=1$. In a way this is the most natural choice of
normalization, because it is possible without any assumptions on $x_0$, $y_0$ beyond $y_0^*x_0=1$. The problem, however, is that it does not define the
eigenvectors uniquely. Suppose we choose $\alpha(\tau)=1/\|x(\tau)\|$, $\beta(\tau)=\|x(\tau)\|$. Then $\alpha$, $\beta^*$,
$\hat x$ and $\hat y^*$ are not analytic in $\tau$, but they are differentiable w.r.t.\ to the real and imaginary parts of $\tau$.
However, we could equally well
multiply $\alpha$ and $\beta$ by any unimodular complex number $e^{i\theta}$, so there are infinitely many different choices of
$\alpha$, $\beta$ that are smooth w.r.t.\ the real and imaginary parts of $\tau$ (though not analytic in $\tau$).
A variant on this normalization is known as real-positive (RP) compatibility \cite{GugOve11}:
it requires $y(\tau)^*x(\tau)$ to be real and positive with $\|x(\tau)\|=\|y(\tau)\|=1$.

 \end{enumerate}

  More generally we could define the normalization in terms of two functions 
  $\psi(\hat x(\tau),\hat y(\tau))=1$ and $\omega(\hat x(\tau),\hat y(\tau)))=1$.  
  Depending on what choice is made, there are several
  possible outcomes:  there could be unique $\alpha(\tau)$ and $\beta(\tau)$ that satisfy the two normalization equations,
  as in cases 1 and 2 above when $e_1^{\T} x_0\ne 0$ and $e_1^{\T} y_0\ne 0$; 
  there could be two choices as in case 3 when $x_0^{\T}x_0\ne 0$; there could be an
  infinite number of choices as in case 4; or there could be no $\alpha(\tau)$ and $\beta(\tau)$ 
  that satisfy the normalization equations.  

Perturbation theory for normalized eigenvectors, using several of the normalizations given above, was
extensively studied by Meyer and Stewart \cite{MeySte88} and by 
Bernasconi, Choirat and Seri \cite{BerChoSer11}.

\subsection{Computation}

Suppose we wish to verify the results of Theorems~\ref{evalpert}, \ref{evecpert} and \ref{eprojpert} computationally --- always a good idea!
Let us consider how to do this in \matlab, where eigenvalues as well as right and left eigenvectors
can be conveniently computed by the function {\tt eig}. 
Let Assumption~\ref{assumption}
hold, assuming for convenience that $\tau_0=0$ and $A(\tau) = A_0 + \tau \Delta A$ for some given matrices $A_0$ and  $\Delta A$ with
$\|A_0\|=\|\Delta A\|=1$.  Take $|\tau|$ sufficiently small that 
only one computed eigenvalue $\tilde\lambda$ of $\tilde A \equiv A(\tau)$ is close to the eigenvalue $\lambda_0$ of $A_0$.
Then, assuming the eigenvalue is not too badly conditioned,
we can easily verify the eigenvalue perturbation formula
\eqref{lamderiv} in Theorem~\ref{evalpert} by computing the finite difference quotient $(\tilde\lambda - \lambda_0)/\tau$
and comparing it with \eqref{lamderiv}.  Since the exact eigenvalue $\lambda(\tau)$ of $A+\tau \Delta A$ is analytic in $\tau$, 
it is clear that, mathematically, the difference between the difference quotient and the derivative should be $O(|\tau|)$
as $\tau \to 0$, but numerically, if $|\tau|$ is too small, rounding error dominates the 
computation instead \cite[Ch.~11]{Ove01}.

Similarly, we can compute right and left eigenvectors $\tilde x$ and $\tilde y$ corresponding to $\tilde\lambda$, normalize these so that
$\tilde y^* \tilde x=1$, compute the eigenprojector $\widetilde \Pi = \tilde x\tilde y^*$, and verify that the matrix difference quotient 
$(\widetilde\Pi - \Pi_0)/\tau$
approximates formula \eqref{eigenprojector} for the eigenprojector derivative given by Theorem~\ref{eprojpert}.

What about the eigenvector derivative formulas? Implementing ``normalization 0" defined in \eqref{xytaudef}, 
we can compute normalizations of $\tilde x$ and $\tilde y$ by
\begin{align}
         \tilde{\tilde x}& = (y_0^* \widetilde\Pi x_0)^{-1/2} \widetilde \Pi  x_0 = \left( \frac{\tilde y^*x_0}{y_0^*\tilde x}\right )^{1/2}\tilde x \label{xtildetilde}\\ 
         \tilde{\tilde y} &= (x_0 ^*\widetilde\Pi^*y_0)^{-1/2} \widetilde\Pi^* y_0  = \left( \frac{\tilde x^*y_0}{x_0^*\tilde y}\right )^{1/2}\tilde y. \label{ytildetilde}
\end{align}
Note that the formulas on the right-hand side of \eqref{xtildetilde}, \eqref{ytildetilde} avoid computing the eigenprojector, which can be advantageous if $n$ is large.
Then we can verify that the difference quotients $(\tilde{\tilde x} - x_0)/\tau$ and $(\tilde{\tilde y}^* - y_0^*)/\tau$ approximate the
eigenvector derivatives \eqref{xderiv} and \eqref{yderiv} given by Theorem~\ref{evecpert}. 

Now consider normalizations 1 to 4 as defined in \S~\ref{subsec:normalize}.
For normalization 1 (respectively, normalization 2) the formulas for the normalized eigenvectors and their derivatives
given by  \eqref{xhatyhat} and \eqref{xhatyhat_derivs} together
with \eqref{alphatau1}, \eqref{betatau1} (respectively \eqref{alphatau1}, \eqref{betatau2})
can be verified easily provided the first
components of $x_0$ and $y_0$ are not zero.
Of course, the index 1 is arbitrary. A better choice is to use $e_j^{\T}\hat x(\tau)=e_k^{\T}\hat y(\tau)=1$, where
$j$ (respectively $k$) is the index of an entry of $x_0$ (respectively $y_0$) with maximum modulus, but this requires
access to $x_0$ and $y_0$.   
In the case of normalization 3, the formulas for the eigenvectors and their derivatives
given by \eqref{xhatyhat}  and \eqref{xhatyhat_derivs} together with \eqref{alphatau3} and \eqref{betatau3}
can also be verified easily, with a caveat due to the freedom in the choice of sign. Specifically,
when $\hat x(\tau)$ and $\hat y(\tau)$ are obtained from the computed vectors $\tilde x$ and $\tilde y$,
it is important to ensure that the signs of $\hat x(\tau)$ and $\hat x(\tau_0)$ (and therefore $\hat y(\tau)$ and $\hat y(\tau_0)$) are consistent; 
this can be done by choosing the signs of the real parts (or the imaginary parts) 
of $e_j^{\T}\hat x(\tau)$ and $e_j^{\T}\hat x(\tau_0)$ to be the same, where $j$ is the index of an entry of $\hat x(\tau_0)$ with maximum modulus.

As for normalization 4, although we could arbitrarily choose $\alpha$ and $\beta^*$ to be smooth functions w.r.t.\ the real and imaginary
parts of $\tau$, there is no way
to know how to obtain smoothly varying computed eigenvectors from the unnormalized computed eigenvectors $\tilde x$ and $\tilde y$.

Summarizing, the formulas for the derivatives of the eigenvalue and the eigenprojector, given in Theorems 1 and 3, are easily verified
numerically, while for eigenvector normalization 0 given by \eqref{xytaudef} and normalizations 1, 2 and 3 defined in \S~\ref{subsec:normalize}, 
the formulas for the eigenvector derivatives can also usually be verified computationally.
However, perhaps surprisingly, there is no panacea when it comes to the eigenvectors.  Normalization 0 always requires access to the
eigenvectors $x_0$ and $y_0$, while the only way to ensure that normalizations 1 and 2 are
well defined is by providing access to $x_0$ and $y_0$ so that indices $j$ and $k$ can be
used instead of index 1 if necessary.  As for normalization 3, it may not be well defined, and even it it is, 
care must be taken to avoid inconsistent sign choices. Finally, normalization 4 is simply not well defined.

Verification of the eigenvalue, eigenprojector
and eigenvector formulas is illustrated by 
publicly available \matlab\ programs.\footnote{https://cs.nyu.edu/overton/papers/eigvecpert-mfiles/eigvecvary\_demo.zip. The main routine is {\tt eigValProjVecVaryDemo.m}}

\section{Multiple eigenvalues} \label{multiple}

Theorems~\ref{evalpert}, \ref{evecpert} and \ref{eprojpert} do not generally hold when $\lambda_0$ is a multiple eigenvalue.
In this section we consider two illuminating examples where multiple eigenvalues enter the picture. 
Recall that algebraic and geometric multiplicity, along with the terms semisimple (nondefective) and nonderogatory,
were defined at the end of \S \ref{evalpert}.
In the following, we use the principal branch of the square root function for definiteness, but any branch would suffice.

\noindent
{\bf Example 1.} 
Let 
\[ 
      A( \tau ) = \left[ \begin{array}{cc} 0 & 1 \\ \tau & 0 \end{array} \right],
 \]
 with $\tau\to \tau_0 = 0$. The limit matrix $A_0=A(0)$ is a Jordan block, with 0  a defective, nonderogatory eigenvalue with algebraic multiplicity 2 and
 geometric multiplicity 1. The corresponding right and left eigenvectors are $x_0=[1,~0]^{\T}$ and $y_0=[0,~1]^{\T}$, but these are mutually orthogonal
 and cannot be scaled so that $y_0^*x_0=1$. For $\tau\not = 0$, $A( \tau )$ has two simple eigenvalues $\lambda_{1,2}(\tau)=\pm\tau^{1/2}$, which are 
not analytic in any neighborhood of 0. The corresponding right and left eigenvectors are uniquely defined, up to scalings, by
\begin{equation} \label{rightleftevecs}
x_{1,2}(\tau) = \begin{bmatrix} 1 \\ \pm \tau^{1/2} \end{bmatrix}~~\mbox{and}~~
y_{1,2}(\tau)^* = \begin{bmatrix} \pm \tau^{1/2} & 1 \end{bmatrix}  .
\end{equation}
The right eigenvectors $x_{1,2}(\tau)$ are not analytic near 0, and they \emph{both} converge to the unique right eigenvector $x_0$ of $A_0$ as $\tau\to 0$.
Likewise $y_{1,2}(\tau)^*$ are not analytic, and they both converge to $y_0^*$
as $\tau\to 0$.
For $\tau\ne 0$, we can scale the eigenvectors so that $y_1(\tau)^*x_1(\tau)=y_2(\tau)^*x_2(\tau)=1$,
but then either $x_j(\tau)$ or $y_j(\tau)$ (or both) must diverge as $\tau\to 0$, for $j=1,2$.

This example is easily extended to the $n\times n$ case, where $A(\tau)$ is zero except for a single superdiagonal of 1's and a bottom left
entry $\tau$, so that $A_0$ is a single Jordan block with 0 a nonderogatory eigenvalue with algebraic multiplicity $n$. The eigenvalues of $A(\tau)$ are
then the $n$th roots of unity times $\tau^{1/n}$. In fact,
Lidskii \cite{Lid66E} gave a remarkable general perturbation theory for eigenvalues of a linear family $A(\tau)$ for which an eigenvalue of
$A_0$ may have \emph{any} algebraic and geometric multiplicities and indeed any Jordan block structure; see also \cite[Sec.~7.4]{Bau85}
and \cite{MorBurOve97}. These results are not described in Kato's books.
However, Kato does treat eigenvalue perturbation in detail in the case that the eigenvalues of $A_0$ are semisimple \cite[Sec.~II.2.3]{Kat82};
see also \cite{LanMarZho03}.
Even in this case, the behavior can be unexpectedly complex, as the following example shows.

\noindent
{\bf Example 2.} \cite[p.~394]{LanTis85}
Let 
\[
    A( \tau ) = \left[ \begin{array}{cc} 0 & \tau \\ \tau^2 & 0 \end{array} \right]
\]
with $\tau\to 0$. This time the limit matrix $A_0=A(0)$ is the zero matrix, with 0 a semisimple eigenvalue with algebraic and geometric multiplicity 2,
so we can take \emph{any} vectors in $\C^2$ as right or left eigenvectors of $A_0$. 
The eigenvalues of $A(\tau)$ are $\lambda_{1,2}(\tau) = \pm \tau^{3/2}$, which are not analytic. The corresponding right and left eigenvectors 
are uniquely defined, up to scaling, by the same formulas \eqref{rightleftevecs} as in the previous example.
So again $x_{1,2}(\tau)$ (respectively $y_{1,2}(\tau)^*$) are not analytic, and both converge to the same vector $x_0$ (respectively $y_0^*$) as
in the previous example when $\tau \rightarrow 0$, although in this case $A_0$
has two linearly independent right (and left) eigenvectors.  There is no right eigenvector of $A( \tau )$
that converges to a vector that is linearly independent of $x_0$ as $\tau\to 0$, and likewise no left 
eigenvector that converges to a vector linearly independent of $y_0$.

In the examples above, a multiple eigenvalue splits apart under perturbation. To avoid dealing with this complexity,
one may study the average behavior of a cluster of eigenvalues, and the corresponding invariant subspace, under perturbation.
There is a large body of work on this topic: see the books by Kato \cite[Sec.~II.2.1]{Kat82}, Gohberg, Lancaster and Rodman \cite{GohLanRod06},
Stewart \cite[Ch.~4]{Ste01} and Stewart and Sun \cite[Ch.~5]{SteSun90}, the unpublished technical report
by Sun \cite[Sec.~2.3]{Sun98}, two surveys on Stewart's many contributions by 
Ipsen \cite{Ips10} and Demmel \cite{Dem10}, and papers such as \cite{BinDemFri08,Dem87,DieFri01,KarKre14}.

In the case of a Hermitian family $A(\tau)$, the perturbation theory for multiple eigenvalues simplifies greatly;
the pioneering results of Rellich were already mentioned in \S~\ref{intro}.  See \cite[Sec.~11.7]{LanTis85} and \cite{GohLanRod85} for more details.

\section{Concluding Remarks}
In this paper we have presented two detailed yet accessible proofs of  first-order perturbation results for a simple
eigenvalue of a matrix and its associated right and left eigenvectors and eigenprojector. We hope this will facilitate the dissemination of these important results
to a much broader community of researchers and students than has hitherto been the case. We have also
tried to convey the breadth and depth of work in the two principal relevant research streams, Analytic Perturbation Theory
and Numerical Linear Algebra. There are, of course, many generalizations that we have not even begun to explore in
this article. Just as one example, nonlinear eigenvalue problems, where one replaces $A-\lambda I_n$ by a matrix function $F(A,\lambda)$ 
with polynomial or more general nonlinear dependence on $\lambda$, arise in many important applications. Perturbation theory for nonlinear eigenvalue problems,
representing the APT and NLA communities respectively, may be found in \cite{AndChuLan93} and \cite{BinHoo13}.
Finally we remark that our bibliography, while fairly extensive, is in no way intended to be comprehensive.

\bibliographystyle{alpha}
\bibliography{EigVecPert_MathSciNet,EigVecPert_Other}

\end{document}